\documentclass[12pt,reqno,a4paper]{article} 
\usepackage[paper=a4paper,left=25mm,right=25mm,top=35mm,bottom=35mm]{geometry}
\usepackage[breaklinks=true,colorlinks=true,linkcolor=black,urlcolor=black,citecolor=black]{hyperref}
\usepackage{lmodern}
\usepackage[T1]{fontenc}
\usepackage[english,ngerman]{babel}
\usepackage[utf8]{inputenc}
\usepackage{indentfirst}
\usepackage{cite}
\makeatletter
\renewcommand\@makefntext[1]{\leftskip=2em\hskip-0.5em\@makefnmark#1}
\makeatother
\usepackage{graphicx} 
\usepackage{float}
\usepackage{amsmath}
\usepackage{amssymb}
\usepackage{mathtools}
\let\oldbibliography\thebibliography
\renewcommand{\thebibliography}[1]{\oldbibliography{#1}
\setlength{\itemsep}{0\baselineskip}}
\begin{document}
\noindent
\begin{center}
	\begin{LARGE}
		\textbf{Günter Hellwig (1926--2004) -- in memoriam}\footnote{This is an extended version of a talk which was given at the memorial colloquium at Aachen on 12 July, 2005. It was privately circulated but not published. The photo was taken in Munich in 2001 at the celebration of Ernst
Wienholtz' 70th birthday.}
	\par
	\end{LARGE}
\end{center}

\begin{center}
	\begin{large}
		\textbf{Hubert Kalf}
	\end{large}\\[1em]
	Mathematisches Institut der Ludwig-Maximilians-Universität München,\\ 
	Theresienstraße 39, 80333 Munich, Germany\\[1em]
	Email: \texttt{kalf@mathematik.uni-muenchen.de}
\end{center}

\begin{figure}[H]
\centering
\includegraphics[scale=0.99]{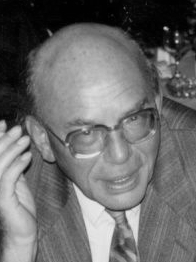}
\end{figure}

Totally unexpectedly Günter Hellwig died on 17 June, 2004. He was the author of influential textbooks on partial differential equations (PDEs) and differential operators of mathematical physics, an enthusiastic and inspiring teacher to generations of engineers until his retirement in 1991, organiser of PDE conferences at Oberwolfach for almost a quarter of a century and a pioneer in index theory.

\paragraph{The War Years and Göttingen.} Born on 9 February, 1926 at Oberschöna near Freiberg, Germany, he belonged to the generation of youths who had to help at anti-aircraft batteries. Drafted in 1943, almost two years before properly finishing his grammar school education, he was engaged for a time at the air-raid shelter near the Zoological Garden of Berlin. It was there that Claus Müller (1920--2008), who was part of a naval research group under the commandership of Helmut Hasse (1898--1979), frequently sought protection. (For Hasse see \cite{Ed,Se}.) It was not until after the war at Göttingen that Müller and Hellwig first met and again much later that they became life-long colleagues at the Rheinisch-Westfälische Hochschule Aachen. 

In 1944 when Hellwig was supposed to be transferred to a parachute unit he had the good fortune of meeting a young man who, because of misbehaviour, had just been dismissed from a training camp named ``Prinz Eugen''. The aim of this camp which was located at the remote place of Stegskopf in the Westerwald was to prepare young people for the study of high frequency techniques. Hellwig contacted this camp and was accepted, which presumably saved his life, because the parachute unit was later wiped out. As he had no higher education entrance qualification (German ``Abitur''), he went to the destroyed city of Dresden early in 1945 and there managed to obtain some form of ``Ersatz'' school certificate which he enjoyed playing on in later life, claiming that he had no proper German ``Abitur''.

Just before the end of the war he was taken prisoner by the US army and brought to a camp near Bad Kreuznach. In order to cope with several hundred thousand prisoners about 18 camps were set up, most of them on meadows along the left bank of the river Rhine between the cities of Bonn and Mainz and along the river Nahe. Food was sparse and the hygienic conditions shocking; thousands died of hunger or dysentery. Owing to the fact that he had some knowledge of English Hellwig was recruited for work in camp administration. He seized the opportunity issuing his own certificate of discharge and having it signed by an officer.

The University of Göttingen reopened in autumn 1945 offering again courses in all of its subjects. The course in analysis for mathematicians and physicists was held by Rellich (1906--1955). (For Rellich see \cite{Cou,Ka}.) The audience was large, almost 300 students. The teaching assistant of the problem session Hellwig attended was Claus Müller. However, Müller left Göttingen for Bonn after only one semester and later went on to Cambridge on a scholarship. He became professor at Aachen in 1955. Hellwig and Heinz-Otto Cordes (1925--2018) competed in solving the starred exercises. Later Cordes became Rellich's assistant; he left Germany in the 1950s and became professor at Berkeley. Another student in Rellich's course to become famous was Martin Beckmann (1924--2017), later Professor of Economics at Brown University and at the same time Professor of Applied Mathematics (Operations Research) at the Technical University of Munich.

In a seminar Hellwig caught Rellich's attention by giving an elegant proof of Lyapunov's inequality \cite[p.~406]{Li}
\begin{align}\label{(1)}
	\int_{a}^{b}\frac{\vert f^{\prime\prime}\vert}{f}>\frac{4}{b-a}\,,
\end{align}
which holds if $f\in C^{0}([a,b])\cap C^{2}((a,b))$  is a function which is positive in the interior and zero at the endpoints of the interval $[a,b]$. (The left-hand side of~\eqref{(1)} may be infinite.) To prove~\eqref{(1)} Hellwig observed that there is a point $x_{1}\in (a,b)$ with $f(x_{1})=\mathrm{max}\,f>0$. The mean-value theorem then
provides numbers $\xi_{1}\in (a,x_{1})$ and $\xi_{2}\in (x_{1},b)$ with 
\begin{align*}
	f(x_{1})=(x_{1}-a)f^{\prime}(\xi_{1})\,,\qquad -f(x_{1})=(b-x_{1})f^{\prime}(\xi_{2})\,.
\end{align*}
Hence
\begin{align}\label{(2)}
	\int_{a}^{b}\frac{\vert f^{\prime\prime}\vert}{f}&>\frac{1}{f(x_{1})}\int_{a}^{b}\vert f^{\prime\prime}\vert\\&\geq \frac{1}{f(x_{1})}\int_{\xi_{1}}^{\xi_{2}}\vert f^{\prime\prime}\vert\nonumber \\ &\geq \frac{1}{f(x_{1})}\bigg\vert\int_{\xi_{1}}^{\xi_{2}} f^{\prime\prime}\bigg\vert  =\frac{b-a}{(x_{1}-a)(b-x_{1})}\geq \frac{b-a}{\frac{b-a}{2}\cdot\frac{b-a}{2}}\nonumber\, .
\end{align}
The inequality in~\eqref{(2)} is strict, otherwise $f$ would be a constant and therefore zero. (The example of a particular concave function which is linear near $a$ and $b$ shows that the constant 4 in~\eqref{(1)} is best possible.) Inequality~\eqref{(1)} can be interpreted as a necessary condition for $-u^{\prime\prime} + qu = 0$ to have a solution with two consecutive zeros. The argument above already shows characteristic features of Hellwig's later work, i.e.~to combine simple tools in very clever ways. (Cordes generalised~\eqref{(1)} in \cite{Cor}.)

Speaking of Rellich's seminars -- Hellwig told the following story: Being asked what else one should sign up for, Rellich replied: ``You will be sufficiently engaged for the whole semester with my seminar.'' 

The mathematician who, beside Rellich, significantly influenced Hellwig was Wolfgang Haack (1902--1994) \cite{Br}. During the war Haack had worked in ballistics and found a shape for a projectile with minimal air resistance in the supersonic region. Theodore von Kármán (1881--1963) had believed that this was impossible to achieve on the basis of the linearised equations of gas dynamics. Immediately after the war Haack acted as advisor to the British Task Force and to the Swiss company Oerlikon. For the summer semester of 1949 he replaced Wilhelm Magnus (1907--1990) at Göttingen who had gone to the United States (it turned out to be permanently) and lectured on PDEs with Hellwig as his teaching assistant. A distinguished student on this course was Jürgen Moser (1928--1999) who had come to Göttingen in 1947 \cite{Ma}.

\paragraph{Berlin.} In autumn 1949 Haack went to the Technical University of Berlin as successor of Georg Hamel (1877--1954) and Hellwig, who had just completed his diploma thesis, joined him as his assistant. Moving to Berlin at that time was not entirely without risk; the Airlift had ended successfully but the future of the city was quite uncertain. In Berlin Hellwig started writing papers in quick succession, habilitating in 1952, one year after his dissertation, at the early age of 26. This was all the more remarkable since the administrative work he had to do for Haack was enormous. On one occasion Hellwig exclaimed ``My assistants are going to have a better life!'' and this was, indeed, to be the case.

Papers \cite{1,2} form the essence of Hellwig's dissertation. By means of a systematic use of Pfaffian forms an analogue of Riemann's integration method was found for first-order hyperbolic systems of linear or quasi-linear equations in two variables. This completed a development which started with Holmgren around 1910. (These results together with the finishing remarks in \cite{14} were also incorporated into the books \cite{16,19,HW}.)

When turning to elliptic systems in the plane Hellwig stumbled upon the following theorem \cite[p.~219]{Hi2} by Hilbert. $C$ is a positively oriented smooth Jordan curve with $J$ as its interior, and the functions $p,q,k$ and $l$ have smooth extensions on $\overline{J}$.  (The translation which follows is mine.) 

\noindent ``When the partial differential equations
\begin{align}\label{(3)}
	\frac{\partial u}{\partial x}-\frac{\partial v}{\partial y}=pu+qv\,,\qquad \frac{\partial u}{\partial y}+\frac{\partial v}{\partial x}=ku+lv\,,
\end{align}
apart from $u=0$, $v=0$, have no system of solutions $u$, $v$ such that $u$ vanishes on the boundary curve $C$, then they have a solution system with the property that $u$ assumes prescribed values $f(s)$ on the boundary curve $C$. Otherwise, i.e.~if there exists a solution system $u$, $v$ of~\eqref{(3)} with the property that $u$ vanishes on $C$ and $u$, $v$ are not both identically zero in $J$, then a solution system $u$, $v$, with $u$ taking prescribed boundary values $f(s)$ on $C$, certainly exists when $f(s)$ satisfies a finite number of certain linear integral conditions.''

Hellwig showed that the first part of the alternative can never occur, which makes Hilbert's theorem vacuous \cite{4}; in fact, the following holds: If $\kappa$ is the winding number (concept and name date back to a paper by Kronecker published in 1869) of the vector field which describes the boundary conditions and $T$ is the linear operator generated by~\eqref{(3)}, then
\begin{align*}
	\mathrm{dim}(\mathrm{ker}\,T)=\begin{cases} 2\kappa+1 & \text{if}\quad\kappa\geq 0\\ 0 & \text{if}\quad\kappa<0\end{cases}\,,\qquad \mathrm{dim}(\mathrm{ran}\,T)=\begin{cases} 0 & \text{if}\quad\kappa\geq 0\\ -(2\kappa+1) & \text{if}\quad\kappa<0\end{cases} \,, 
\end{align*}
i.e., $\mathrm{ind}\,T=\mathrm{dim}(\mathrm{ker}\,T)-\mathrm{dim}(\mathrm{ran}\,T)=2\kappa+1$. This was proved by Hellwig for $\kappa=0$ \cite{5}. The general situation can be reduced to this special case which was effected in \cite{Ha2}. In the context of singular integral operators this index had been introduced by Fritz Noether (1884--1941) in \cite{No} where
he corrected an erroneous remark made by Hilbert in \cite[p.~239]{Hi1}. The integral operators considered by Fredholm in his famous investigations have the property $\mathrm{ind}\,T=0$.

In 1920 Liebmann (1874--1939) proved that a closed convex surface is rigid and he conjectured that it becomes deformable when an arbitrarily small piece is cut out. Cohn-Vossen (1902--1936) and Süss (1895--1958) independently of each other investigated infinitesimal deformations in papers which appeared in 1927. (For Süss see \cite{Re,Se}.) Süss observed that the fundamental equations describing the surfaces in question are elliptic systems of the form Hilbert had considered and inferred their solvability from Hilbert's remarks above. (In contrast, Cohn-Vossen relied on the solvability
of a scalar elliptic equation. Both proofs are described by Rembs (1890--1964) and Grotemeyer (1927--2007) in \cite[p.~200 f.]{Ef}; the references which are not given here can also be found in this book.) So it was natural to take up the matter again. Meanwhile I.~N.~Vekua's (1907--1977) more general results from the function-theoretic point of view on elliptic systems in the plane had appeared (in contrast to \cite{5,6} and \cite{Ha2} there was no restriction of the size of the underlying domain). With the help of \cite[\S 8.10]{Ve} Hellwig was able to prove Liebmann's conjecture in \cite{11}, at least in the class of surfaces which are four times Hölder continuously differentiable.

Papers \cite{7,8,9} and \cite{12} are concerned with hyperbolic equations of two variables which become parabolic on the boundary. \cite{7} connects such problems with the spectral properties of the Sturm-Liouville operator
\begin{align*}
	L:=\frac{1}{k}\bigg(-\frac{d}{dx}p\frac{d}{dx}+q\bigg)\, .
\end{align*}
The functions $p$, $q$ and $k$ should at least be continuous on some open interval $I$ (which may be unbounded) and $p$ and $k$ are supposed to be positive. Let $r_{1},r_{2},r_{3}\in\mathbb{R}$. Then the equation
\begin{align}\label{(4)}
	Lu+r_{1}u_{yy}+r_{2}u_{y}+r_{3}u=f\qquad \text{on}\qquad I\times (0,\infty)
\end{align}
is hyperbolic if $r_{1}>0$ ($L$  acts on the first variable of $u$).  Suppose $L$ is self-adjoint in the Hilbert space of functions which are square-integrable with weight $k$ and its spectrum is bounded from below. (The first assumption is intimately connected with the behaviour of the characteristic directions of~\eqref{(4)}, given by $(k/p)^{1/2}$, near the boundary of $I$.) Given suitable decay estimates for the resolvent of $L$, Hellwig shows that the initial-value problem for~\eqref{(4)} can be solved by means of the Laplace transform of this resolvent.

The situation that $I$ is bounded and $L$ is a regular Sturm-Liouville operator with coupled homogeneous or inhomogeneous boundary conditions is considered in \cite{13,15} for the parabolic case $r_{1}=r_{3}=0$, $r_{2}>0$. Here the required decay estimates can be provided by means of a careful asymptotic analysis of a fundamental system of solutions.

More general equations of the form
\begin{align}\label{(5)}
	(a(x,y)u_{x})_{x}-(c(x,y)u_{y})_{y}-d(x,y)u=f(x,y)
\end{align}
which are assumed to be hyperbolic for $x>0$ and parabolic for $x=0$ are considered in \cite{9}. Here a generalised solution of the initial-value problem on the parabolic line is defined via a singular Volterra equation, which is solved by iteration. This problem was also the topic of his talk at the International Congress of Mathematicians in Amsterdam in 1954 \cite{10}.

Invited by Courant (1888--1972), he spent the academic year 1954/55 on a Fulbright grant at the Institute of Mathematical Sciences of New York University, the Mecca for researchers in PDEs. Here life-long contacts developed with Bers (1914--1993), Courant, Friedrichs (1901--1982), John (1910--1994), Lax (1926--2025), Morawetz (1923--2017) and Nirenberg. An outgrowth of his stay at New York is his massive paper \cite{12}, where the method of \cite{9} is used to obtain unique (generalised) solutions for five types of initial-value problems when~\eqref{(5)} becomes parabolic on the $y$-axis.

In 1958 he was offered a full professorship at the University of Würzburg, which he declined, but he accepted a call to the newly created fourth chair of mathematics at the Technical University of Berlin in the same year. This new position enabled him to engage as assistant Ernst Wienholtz (1931--2003), Rellich's last doctoral student who after Rellich's untimely death had just been promoted at Göttingen by Carl Ludwig Siegel (1896--1981). In his dissertation Wienholtz proved among other things that a Schrödinger operator which is bounded from below is already essentially self-adjoint, at least if the potential is continuous, answering a question Rellich had posed at the Amsterdam Congress. Only much later was it realised in the West that this result had already been obtained by A.~Ya.~Povzner (1915--2008) in 1953. Wienholtz' simple and elegant proof was, however, taken into Glazman's book \cite[p.~58 f.]{Gl}, and results of this kind are sometimes referred to as Povzner-Wienholtz type theorems.

In 1960, Hellwig had just turned 34, his great book on PDEs appeared. An enlightening brief overview over the classic triad, the wave, Poisson and heat equation, is given before the problem of a division into types is addressed in a more systematic way. Questions of uniqueness of solutions, classical or generalised, follow. (It needs to be noted that a tool like the Lax-Milgram lemma was only found some 5 years before.) A great deal of his own research is incorporated, as are two new results by Wienholtz, the extension of the regularity result of Eberhard Hopf (1902--1983) to weak solutions of a general second-order elliptic equation \cite[IV.4.2]{16}, and his ingenious proof that such solutions assume their boundary values continuously at a smooth portion of the boundary \cite[IV.3.7--8]{16}. This book was an immediate and a lasting success. In the year of his death he told me with some chagrin that he had found out that the then recent and now popular book by Evans \cite{Ev} did not cite his book.

The well-known Haack-Hellwig Oberwolfach Conferences on PDEs started in March 1961. It was the success of this first conference and Courant's eventual enthusiastic reaction which helped to secure the future of the Institute
at Oberwolfach, whose financial support at that time was uncertain. After Haack's retirement Hellwig organised the conferences from 1973--1979 with Erhard Heinz (1924--2017) and then with Joachim Weidmann until 1985 when the tradition was discontinued in a manner which deeply hurt Hellwig.

In 1961 Hellwig married Birgitta Öman (1932--2021), a young Swedish mathematician whom he had first met at the ICM in Edinburgh in 1958. They had a son and three daughters one of whom also became a mathematician working in graph theory.

Working independently, both Rellich and Vekua had shown the following result in 1943: If u is a non-trivial function which satisfies
\begin{align}\label{(6)}
	-\Delta  u=\lambda u
\end{align}
for some $\lambda>0$ in the exterior of a ball, then there exist positive numbers $p$, $s$ such that
\begin{align}\label{(7)}
	\int\limits_{s<\vert x\vert<t}\vert u(x)\vert^{2}\,dx\geq pt\qquad (t>s)\, .
\end{align}
In particular, there can be no non-trivial solution which is of integrable square in an exterior domain. The estimate~\eqref{(7)} is best possible in the sense that there are simple examples where the left-hand side of~\eqref{(7)} grows linearly.

In order to establish the non-existence of positive eigenvalues for a larger class of operators, it was important to dispense with the tools Rellich had used in his proof of~\eqref{(7)}, viz.~the expansion of $u$ into spherical harmonics and the asymptotic behaviour of Bessel functions. In \cite{17} Hellwig gave an elegant proof of~\eqref{(7)} which used the mean-value property of solutions of~\eqref{(6)} only, but he had to assume that they are defined on the whole of $\mathbb{R}^{n}$. This proof was also inserted into the English translation of his PDE book \cite[p.~113 f.]{19}, an edition which he particularly liked because of its elegant layout. The extension of~\eqref{(7)} to more general situations required indeed very sophisticated means and is primarily due to Agmon (1922--2025) and Kato (1917--1999). References can be found in \cite{EK}.

Paper \cite{17} indicated already his shift of interest towards problems Rellich had been concerned with. This became even more apparent in his second book \cite{20}, which was indeed dedicated to the memory of Rellich. For quite some time this was the only textbook from which one could learn Weyl's ingenious theory of the limit-point and limit-circle case as well as the description of the self-adjoint extensions of a Sturm-Liouville operator in terms of the Rellich initial-numbers and the distinguished role the Friedrichs extension plays. As far as Schrödinger operators are concerned, the central idea of the book is to define them on a vector space of smooth functions on which they are essentially self-adjoint and where their spectral resolution takes place, i.e.~where their eigenfunctions as well as their eigenpackets lie. His wife helped with this book and she also translated it into English; this edition also contains her generalisation of Levinson's self-adjointness criterion \cite[\S6.2]{21}. This edition was prepared when they were still at Berlin but was only published after their move to Aachen.

\paragraph{Aachen.} In 1966 Hellwig succeeded Hubert Cremer (1897--1983) in the chair of mathematics at RWTH Aachen. His institute was big; there were two secretaries and ten positions for assistants. The number of students to deal with -- originally students of all engineering faculties, somewhat later electrical engineers and physicists -- was huge. At some point 1,200 students had signed up for his four-semester course ``Higher Mathematics''. The largest lecture theatre could seat one thousand people; therefore, quite a few had to put up with sitting on stairs. The two youngest assistants had to accompany him -- looking after the lighting, the microphone and the two overhead projectors, taking notes of the lecture, correct slips of his pen (rarely necessary) and just being there in case anything untoward happened.

During the student unrests after 1968 there were some occasional incidents, but usually disturbances came from the outside. The vast majority of his students liked the lively way in which he presented the material and appreciated the genuine concern he had for them. Those who failed the final written examination twice had to be examined verbally. This was done in groups of ten. The questions Hellwig asked had to be answered in written form. Since assistants were sitting amidst the candidates -- and could also lend a helping hand -- in the end some sort of individual assessment was possible. (Much later this system had to be abandoned for legal reasons.)

It is not hard to imagine how much these routine tasks strained a sensitive person such as Hellwig. It must also be said that we, as his assistants, were frequently not easy to handle either when we felt driven to raise the level of sophistication or rigour which, looking back, I find was remarkably high. To relax, Hellwig went out for long walks on his own enjoying the wide vistas of a highland moor on the Belgian side of the border (Hohes Venn).

The highlight of the week was the Hellwig research seminar on Tuesdays after lunch. Usually someone talked about a result he had  just found himself or read about. Hellwig made critical remarks and insisted on details, so that the talk often continued over several sessions. The appended list of doctoral theses and habilitations gives some idea of the variety of topics that were addressed. Naturally, the success of the seminar was not due to Hellwig alone; while it may be inappropriate to single out participants still alive, I feel I must mention the prominent role played by Johann Walter (1932--2008).

After the seminar or the colloquium of the department, Hellwig seized the opportunity to talk. Vividly and with his strong Saxon accent he loved to talk about mathematics and the many mathematicians he knew. He was a great raconteur -- with a dry sense of humour. On one occasion, for instance, he recounted with relish how his secretary in Berlin, when typing the invitations for the winter ball of the Berlin Mathematical Society, made a typing mistake, putting ``widows' ball'' instead. He did not correct her.

The atmosphere at his institute was wonderful. We could work in total freedom on problems that were close enough to allow stimulating discussions and sufficiently separated to prevent rivalries. The 1970s and 1980s were a ``golden age'' of spectral theory and some contributions, albeit in a small area, came from Hellwig's institute at Aachen.

At the end of the 1960s, students generally started calling for lecture notes. Eventually, Hellwig gave in and published at least a small part of his course as a book. It was not without reason that he regarded a manipulation of indefinite integrals as a manipulation of sets of functions. In contrast to his usual pragmatism he developed an elaborate calculus to deal with indefinite integrals \cite[\S5.2]{22} which, however, did not find the attention he had hoped for. The book has a charming little paragraph where, in a playful way, density results and inequalities are derived from very simple observations \cite[\S6.13]{22}.

From the mid 1970s he regularly gave a course of lectures on the relation between science and religion, which became quite popular. He held the beauty of mathematical descriptions of natural phenomena as indicative of divine providence. He loved to detail this with the example of the wave equation which admits travelling waves of arbitrary shape and precludes reverberating effects precisely in three space dimensions.

On a June day in 2004 he went out for a swim -- alone, as was his habit -- but did not return home. He had succumbed to heart failure.

\paragraph{Acknowledgements.} I gratefully acknowledge the help I received when preparing this obituary: Mrs Hellwig provided me with detailed information about her husband's vicissitudes of life during and after the war; H.~O.~Cordes added reminiscences of the Göttingen and H.-W.~Rohde (1936--2021) of the Berlin years. J.~Bemelmans assisted with the list of dissertations and habilitations, and W.~N.~Everitt (1924--2011) put a critical eye to the first draft and made valuable suggestions.
\begin{otherlanguage}{ngerman}
\renewcommand{\refname}{Publications by Günter Hellwig}

\section*{Dissertations}

\begin{list}{}%
  {\setlength{\leftmargin}{2em}
   \setlength{\itemindent}{-2em}
   \setlength{\itemsep}{0.8ex}}
   
\item \textbf{J.~Witte}, Über die Regularität der Eigenfunktionen und Eigenpakete eines
singulären elliptischen Differentialoperators (8.5.1967).

\item \textbf{K.-H.~Jansen}, Kriterien für das Fehlen quadratisch integrierbarer Lösungen der Dif-ferentialgleichung $-\Delta v = f(x, v)$ in Außengebieten (21.6.1968).

\item \textbf{U.-W.~Schmincke}, Über das Verhalten der Eigenfunktionen eines singulären elliptischen Differentialoperators (9.7.1968).

\item \textbf{D.~Kniepert} (1934--2012), Globale Aussagen über die Eigenfunktionen eines singulären elliptischen Differentialoperators im Gebiet $G$ (16.12.1968).

\item \textbf{M.~Teuffel}, Regularitätsaussagen mit Anwendung auf die Spektraltheorie singulärer elliptischer Differentialgleichungssysteme (7.7.1969).

\item \textbf{R.~Wüst} (1943--2021), Stabilität der Selbstadjungiertheit gegenüber Störungen (14.5.1970). 

\item \textbf{H.~Kalf}, Über die Selbstadjungiertheit halbbeschränkter gewöhnlicher oder elliptischer Differentialoperatoren mit stark singulärem Potential (12.2.1971).

\item \textbf{M.~Pisters}, A priori Abschätzungen und Regularitätsaussagen bei elliptischen Differential-gleichungen 2. Ordnung und deren Anwendung auf die Störungstheorie (2.7.1971).

\item \textbf{H.-W.~Goelden}, Über die Nichtentartung des Grundzustandes bei Schrödinger-Operatoren (23.11.1976).

\item \textbf{H.~Bachmann}, Eine Spektraltheorie des 't Hooft'schen Eigenwertproblems im Hilbert-Raum (16.2.1983)
\end{list}
\section*{Habilitations}

\begin{list}{}%
  {\setlength{\leftmargin}{2em}
   \setlength{\itemindent}{-2em}
   \setlength{\itemsep}{0.8ex}}
  
\item \textbf{E.~Wienholtz}, Das Weylsche Lemma für elliptische Systeme partieller Differentialoperatoren oder zur Regularität schwacher Lösungen für elliptische Systeme partieller Differentialoperatoren (24.4.1963).

\item \textbf{H.-W.~Rohde}, Kriterien zur Selbstadjungiertheit elliptischer Differentialoperatoren mit singulärem Potential (18.10.1967).

\item \textbf{J.~Walter}, Diskretes Spektrum bei Sturm-Liouville-Operatoren (13.1.1971).

\item \textbf{U.-W.~Schmincke}, Ausgezeichnete selbstadjungierte Erweiterungen von Dirac-Operatoren mit stark singulärem (insbesondere Coulomb-) Potential (13.12.1972).

\item \textbf{R.~Wüst}, Ausgezeichnete selbstadjungierte Fortsetzungen von Dirac-Operatoren mit stark singulärem Potential, konstruiert mit cut-off-Methoden und einem Spektrallückensatz (10.12.1975).

\item \textbf{H.~Kalf}, Ein Kriterium für das Fehlen von Eigenwerten bei Schrödingeroperatoren (15.12.1976).

\item \textbf{G.~Dziuk}, Über die Glattheit des freien Randes bei Minimalflächen (7.7.1982).

\item \textbf{G.~Ströhmer}, Variationsprobleme mit Nebenbedingungen und der Satz vom Wall (9.5.1984).
\end{list}
\end{otherlanguage}
\end{document}